\documentclass[12pt]{article}

\usepackage{latexsym}
\usepackage[frenchb]{babel}
\usepackage[latin1]{inputenc}
\usepackage{amsfonts,amsmath,amssymb}
\author{Andrea Surroca}
\title{\Large Siegel's theorem and the $abc$ conjecture}
\date{}

\newtheorem{conj}{\textbf{Conjecture}}

\newtheorem{thm}{\textbf{Theorem}}

\newtheorem{lemma}{\textbf{Lemma}}

\newtheorem{prop}{\textbf{Proposition}}
  
\newtheorem{hyp}{\textbf{Hypothesis}}

\newcommand\rat{\mathbf{Q}}

\newcommand\Qbarre{\overline{\mathbf{Q}} }
  
\newcommand\card{\mathrm{card}}

\newcommand\rad{\mathrm{rad}}

\newcommand\supp{\mathrm{supp}}

\newcommand\mult{\mathrm{mult}}

\newcommand\ord{\mathrm{ord}}

\newcommand\espproj{{\mathrm{I} \hspace{-0,1 cm} \mathrm{P}}}

\newcommand\Spec{\mathrm{Spec}}

\newcommand\Div{\mathrm{Div}}

\newcommand\vp{\mathfrak{p}}

\newcommand\vq{\mathfrak{q}}

\newcommand\esp{\hspace{0,2cm}}

\newcommand\Belyi{Bely\u\i}

\begin{document}

\maketitle

\begin{quote}
\textbf{Abstract.}
{\small Following N. Elkies  \cite{elkies} we show that the $abc$ conjecture of Masser-Oesterlé implies an effective version of Siegel's theorem about integral points on algebraic curves, i.e. an upper bound for the $S$-integral points where the dependence on $S$ is explicit. The converse statement is also announced in this note. For both results, the main geometric tool is a theorem of G.V. \Belyi.
 }

\end{quote}

{\small 2000 Mathematics Subject Classification. Primary: 11G30; Secondary: 11J25, 11G50.}

\section{Introduction.}

One of the fundamental theorems in diophantine geometry, proved by C.L. Siegel, says that an algebraic affine curve of genus $\geq 1$ or of genus zero with 3 points at infinity and defined over a number field $K$ has only finitely many integral points. Thanks to the contribution of K. Mahler, Siegel's theorem was extended to the $S$-integral points, that is, those having the prime factors of the denominators in the finite set of valuations $S$ of $K$.

For curves of genus $\geq 2$, Siegel's theorem is supperseeded by Falting's one (Mordell's conjecture), which asserts that there are finitely many rational points on an algebraic projective curve of genus $\geq2$. However, there are no effective versions of these theorems, in the sense that the proofs do not provide an algorithm to find such points. Some effective results concerning integral points are known, thanks to Baker's method based on linear forms in logarithms. Infact, this method gives upper bounds for the height (i.e. the size) of the ($S$-)integral points.

N. Elkies proposed in \cite{elkies} an effective method to study rational points on algebraic curves, but it remains conjectural, because it is based on the $abc$ conjecture of Masser-Oesterlé. This approach was the motivation for the work that we explain here. 

\begin{conj}{($abc$)}\label{abc}

 Let $K$ be a number field. For every $\varepsilon >0$, there exists a real number $c_{\varepsilon, K}>0$ such that, for every $a, b, c \in K \setminus \{0\}$ and verifying $a+b=c$, we have
$$h_{K}(a:b:c)< (1+\varepsilon) \rad (a:b:c) +c_{\varepsilon, K}.$$
\end{conj}

In this note we will  explain some links between upper bounds on the height of $S$-integral points on algebraic curves and the $abc$ conjecture. 

\bigskip

\textbf{Acknowledgements.}
{\small The results quoted in this note come from the PhD thesis of the author, who wants to thank her advisors, Marc Hindry and Michel Waldschmidt. She wants also to thank Machiel van Frankenhuysen for sending her some of his off prints, and the anonimous referee for pointing out to her the references \cite{song-tucker-dirichlet}, \cite{song-tucker-discriminants} and \cite{tucker}. 
 }

\section{Notations.}

Let us denote by $U$ an algebraic affine curve of genus $g$ defined over a number field $K$, $\tilde{U}$ will denote its normalization and $C$ a projective curve containing  $\tilde{U}$. Let $U_{\infty}= C \setminus \tilde{U} = \{P_{1}, \ldots, P_{t}\}$ be the "points at infinity".  The Euler-Poincaré characteristic of $U$ is the number 
$\chi(U)= 2- 2g -t$. The curve $U$ verifies the hypothesis  of Siegel's theorem if and only if $\chi(U)<0$.

Let $S$ be a finite set of valuations of the number field $K$ and denote $O_{K}$ and $O_{K,S}$ respectively, the ring of integers and the ring of $S$-integers of $K$.

Denote also by $U(O_{K,S})$ the $K$-rational points on $U$ which have $S$-integral  coordinates with respect to  the embedding $\phi_{D} : U \hookrightarrow \mathbf{A}^{n}$ given by a very ample divisor $D \in \Div (C)$, such that $\supp (D) = U_{\infty}$.

All along this paper, $h$ will denote the absolute logarithmic Weil's height on $\mathbf{P}^{1}$, $h_{K} = [K:\rat]\,h$, and $h_{U}$ a height on $U$ relative to a degree 1 divisor (e.g. $h_{U}= \frac{1}{\deg(f)}h \circ f$, where $f$ is a non constant regular function on $U$).

If $x=(x_{1},x_{2},x_{3}) \in \mathbf{P}^{2}(K)$, we define the radical (or the support) of $x$ (with logarithmic notation) by
$$\rad (x) = \sum \log N_{K/\rat}(\vp),$$
where the sum is taken over  the prime ideals $\vp$ of $K$ for which the cardinality of $\{v_{\vp}(x_{0}), v_{\vp}(x_{1}), v_{\vp}(x_{2})\}$ is greater than or equal to 2.

For example, if $x_{1},x_{2},x_{3}$ are nonzero rational coprime numbers, $h(x_{0}:x_{1}: x_{2}) = \log \max \{|x_{0}|, |x_{1}|, |x_{2}| \}$ and $\rad(x_{0}:x_{1}: x_{2}) = \log \prod p$, where the product is taken over the prime numbers $p$ dividing $x_{0}, x_{1}$ or $x_{2}$.

\section{Statement of the results.}

Following the ideas proposed by Elkies in \cite{elkies} (see also \cite{franken}) we show that (an effective version of) the $abc$ conjecture implies an effective version of Siegel's theorem i.e. an upper bound for the $S$-integral points where the dependence on $S$ is explicit. (Cf. the last section for the proof.)

\begin{thm}\label{abc-siegel}

Let $K$ be a number field, $U$ an affine algebraic curve defined over $K$ such that $\chi(U)<0$, $S$ a finite set of absolute values of $K$ and $h_{U,K}$ a height function on $U(K)$ associated to a degree 1 divisor.

Suppose the $abc$ conjecture is true.

For every  $\varepsilon \in ]0,-\chi(U)[$ and every $S$-integral point $x$ on $U$, there exists a number $\gamma$ depending on the number field $K$, on the curve $U$, on the choice of the height $h_{U,K}$ and on $\varepsilon$, but neither on the point $x$ nor on the set $S$, such that 
\begin{equation}
h_{U,K}(x) \leq -\frac{1}{\varepsilon + \chi(U)} \sum_{\vp \in S} \log N_{K/\rat}(\vp) + \gamma. \label{ineg.abc-siegel}
\end{equation}
Moreover, if the constant $c_{\varepsilon,K}$ of the conjecture \ref{abc} were  
effective,  the same would be true for  $\gamma$. 
 
\end{thm}

If we consider a version of the $abc$ conjecture 
where the dependence of $c_{\varepsilon, K}$ on the number field  is explicit, e. g. $c_{\varepsilon, K}=2 \log d_{K}+ \kappa$ (cf. \cite{masser}), we would obtain the same upper bound, but with
\begin{equation}
\gamma = -\frac{2}{\varepsilon + \chi(U)} \log d_{K} + \kappa, \label{gamma}
\end{equation}
where $d_{K}$ is the absolute value of the discriminant of $K$
and $\kappa$ only depends on the degree $[K:\rat]$ of the number field $K$.

Inequality (\ref{ineg.abc-siegel}) modified with (\ref{gamma}) suggests the following hypothesis.

Fix an algebraic affine curve $U$ defined over $K$ and verifying $\chi(U)<0$, 
as well as a height function $h_{U}$ on  $U$, associated to a divisor of degree 1.

\begin{hyp}{{\it Siegel $(U,K)$}}\label{se(U,K)} 

For every $\delta \geq [K:\rat]$, there exists constants $k_1(U,h_{U},\delta) >1, \esp k_2(U,h_{U},\delta) >1 $ and $k_3(U,h_{U},\delta) >0$ such that, for every finite extension $L/K$ of degree $[L:\rat] \leq \delta$, for every finite set  $S$ of valuations of $L$ and for every $S$-integral point $x$ on $U$, we have
\begin{equation}\label{upper}
h_{U,L}(x) \leq k_1 \sum_{\mathfrak{q} \in S}  \log N_{L/\rat}(\mathfrak{q}) + k_2 \log d_L + k_3,
\end{equation}
where $d_L$ is the absolute value of the discriminnant of the field $L$ and $h_{U,L} = [L:\rat]\, h_U$.  
\end{hyp}

A similar hypothetical upper bound for rational points was considered by L. Moret-Bailly \cite{moret-bailly} to show that an effective version of Mordell's conjecture would imply the $abc$ conjecture. Moreover, the analogue of the hypothesis \ref{se(U,K)} in the function field case has already been proved (see the appendix of \cite{mathese} and the references therein). 

We prove the following theorem.

\begin{thm}\label{siegel-abc}

Let $K$ be a number field, $U$ an algebraic affine curve defined over $K$ and such that $\chi (U) < 0$.  Let $h_U$ be a height function on $U$ associated to a degree 1 divisor. 

Assume the hypothesis \ref{se(U,K)}.

Then, there exists real positive numbers $\eta_{1}$ and $\eta_{3}$ depending on $U$ and on the degree $[K:\rat]$, and $\eta_{2}$ depending only on $U$ such that, for every $a, b, c \in K \setminus \{0\}$ and verifying $a+b=c$, we have
$$h_{K}(a:b:c)< \eta_{1} \rad(a:b:c) + \eta_{2} \log d_{K} + \eta_{3}.$$

\end{thm}

The proofs of theorem \ref{abc-siegel} and theorem \ref{siegel-abc} make use of a \Belyi \esp function.

\begin{thm}{(\Belyi)}\label{thm.belyi}

An algebraic projective curve $C$ is defined over $\Qbarre$ if and only if there exists a finite and surjective morphism $f : C \rightarrow \mathbf{P}^{1}$, unramified outside $\{0,1, \infty\}$. 

Moreover, if $\Sigma \subset C(\Qbarre)$, we can choose $f$ such that $f(\Sigma)\subset \{0,1,\infty\}$.

\end{thm}

 The construction of such a function is completely explicit and if the curve $C$ is defined over a number field $K$ the same will be true for the function $f$. 
 
\smallskip

For the proof of theorem \ref{siegel-abc} the main idea is to put $y= (a:c) \in \mathbf{P}^{1}(K)$ and prove that a point $x$ in $U$ such that $f(x)=y$, where $f$ is a \Belyi \esp function, is in fact $S'$-integral for some finite extension $L/K$ and some finite set $S'$ of valuations of $L$. To conclude we apply the hypothesis \ref{se(U,K)} to $x$. (Cf. \cite{mathese} for a detailed proof.)

\smallskip

Yu. Bilu obtained in \cite{bilu} an upper bound weaker than (\ref{upper}), but inconditional, for Galois coverings. His result allows us to prove the following theorem.

\begin{thm}\label{abc.exp}

For every number field $K$, there exists real and effectively computable numbers $\gamma_{1}, \gamma_{2}>1$, depending on the degree $[K:\rat]$ and on the discriminant $d_K$ of $K$, such that, for every  $(a,b,c)$ in $K \setminus \{0\}$ and verifying $a+b=c$, 
$$h_K(a:b:c) \leq \exp \{\gamma_{1} \, \rad(a:b:c)  + \gamma_{2}\}.$$

\end{thm}

Theorem \ref{abc.exp} generalizes to any number field a result of Stewart-Yu (over $\rat$) \cite{stewart-yu} which is the best known result in the direction of the $abc$ conjecture.

\section{Proof of theorem \ref{abc-siegel}.}

Let $U$ be an algebraic affine curve  of genus $g$ defined over a number field $K$,  $\tilde{U}$ its normalization and $C$ a projective curve containing $\tilde{U}$. 
Suppose that $\chi(U) <0$. 

\begin{lemma}\label{lemma.belyi}

There exists a rational function $f \in K(C)$ such that 
$$\card(\{x\in U(\Qbarre) / f(x)\in \{0, 1, \infty\}  \}) = \chi(U) + d,$$
where $d$ is the degree of $f$.

\end{lemma}

\noindent \textbf{Proof of the lemma \ref{lemma.belyi}.}

Apply theorem \ref{thm.belyi} to the projective curve $C$, and to the set $\Sigma = \{P_{1}, \ldots, P_{t}\} $. Fix such a \Belyi \esp function $f$ and set $d = \deg(f)$.

Set also $m = \card(\{x\in C(\Qbarre) / f(x)\in \{0, 1, \infty\}  \})$.

The function $f$ is ramified only above $0,1$ and $\infty$, so the Riemann-Hurwitz formula gives
$$\sum_{P\in C, f(P)\in \{0,1,\infty\}} (e_{P}(f)-1)  = 2d + 2g - 2,$$
 and then 
 $$m= 3d -\sum_{P\in C, f(P)\in \{0,1,\infty\}} (e_{P}(f)-1)  = d-2g+2. $$

Moreover, $f(\{P_{1}, \ldots, P_{t}\}) \subset \{0,1,\infty\}$, then 
$$\card(\{x\in U(\Qbarre) / f(x)\in \{0, 1, \infty\}  \}) = m - \card(U_{\infty}(\Qbarre)) = m - t =  d +\chi(U).$$
\hfill $\Box$
 
Let $S$ be a finite set of absolute values of $K$.
Let $x$ be in $U(O_{K,S})$ such that $f(x)=(r:1)\in \mathbf{P}^{1}\setminus \{0, 1, \infty\}$. (If $f(x)\in \{0,1,\infty\}$, then $h_{K}(f(x))=0$ and the inequality (\ref{ineg.abc-siegel}) is verified.) We can then apply conjecture \ref{abc} to $P_{x}=(r:1-r:1)$. For every $\varepsilon >0$, there exists a constant $c_{\varepsilon, K}$, such that
$$ h_{K}(P_{x}) \leq (1+\varepsilon)\rad (P_{x}) + c_{\varepsilon, K}.$$
Choose $h_{U,K}= \frac{1}{d}(h_{K}\circ f)$.  We then have 
$$d \, h_{U,K}(x) =h_{K}(f(x)) = h_{K}(r:1) \leq h_{K}(r:1-r:1) = h_{K}(P_{x}).$$

We conclude with the following proposition which gives an upper bound of the radical in terms of the height.

\begin{prop}\label{radical} 

There exists constants $c_1$ and $c_2$ depending only on the degree $[K:\rat]$ of the number field $K$ and of the curve $U$, such that 
$$\rad(P_{x}) \leq \left( \frac{\chi(U)}{d} + 1 \right) h_K(f(x)) + c_1 \sqrt{h_K(f(x))} + c_2 + \sum_{\vp \in S} \log N(\vp). $$        
\end{prop}

To prove proposition \ref{radical}, let us introduce, like in \cite{franken}, the divisors of zero, of 1 and of the point at infinity, with respect to the \Belyi's function $f$. Let $D_{0}=f^{\ast}(0) = \sum_{x\in C(\Qbarre), f(x)=0} e_{x}(f)(x) \in \Div(C)$, where   $e_{x}(f)$ is the ramification index of $f$ at the point $x$, be the divisor of zero. Denote by $D_{0}|_U = \sum_{x\in U(\Qbarre), f(x)=0} e_{x}(f)(x)$ the restriction of the divisor to the affine curve $U$. In the same way, define the divisor $D_{1}$ of 1, and $D_{\infty}$  of the infinity. 
The curve $U$ is defined over $K$, so \Belyi's theorem allows us to take $f$ defined over $K$. The divisors $D_{0}|_{U}, D_{1}|_{U}$ and $D_{\infty}|_{U}$ are defined over $K$ and have a decomposition into irreducible divisors:
$$D_{0}|_U = e_{1}M_{1} + \cdots + e_{i}M_{i}, \esp \esp \esp$$
 $$D_{1}|_U = e_{i+1}M_{i+1} + \cdots + e_{j}M_{j},$$
 $$D_{\infty}|_U = e_{j+1}M_{j+1} + \cdots + e_{k}M_{k}.$$

 By the Riemann-Roch theorem we can prove that for every $\nu \in \{1,\ldots, k\}$, there is a function  $m_{\nu} : C \longrightarrow \espproj^1$ such that 
 $$(2g+1) M_{\nu}= m_{\nu}^{\ast}(0).$$

We can show that there exists a finite set  $T$ of valuations of $K$ (depending only on the curve $U$) such that,  there exists a model $\mathcal{C}$ of the curve $C$ over $\mathcal{B}$, where  $\mathcal{B} = \Spec (O_{K,T}) = \Spec(O_K) \setminus T$, such that:\\
\indent i) all the special fibers of this model are smooth projectives curves,\\
\indent ii) the functions $f$ and $m_{\nu}$ extend to morphisms  $\tilde{f}$ and $\tilde{m_{\nu}}$ from $\mathcal{C}$ to $\espproj^1_{O_{K,T}} = \espproj^1_\mathcal{B}$ over $\mathcal{B}$. 

In other words, we can say that the curve $C$ has good reduction at primes outside $T$, and that the functions $f$ and $m_{\nu}$ are not reduced to the identically zero function.

\begin{lemma}\label{log}
 
Let $\vp$ be a prime  ideal in $\mathcal{B}$, which contributes to the radical of  $P_x$. 
If $\vp\notin S$, then 
 $$\log N(\vp) \leq \frac{d_{\vp}}{2g+1} \sum_{\nu =1}^k \max \{0, v_{\vp}(m_{\nu}(x)) \}.$$
\end{lemma} 

\noindent \textbf{Proof of the lemma \ref{log}.}

Fix $\vp\in \mathcal{B}$ and put $\mathcal{U} = \mathcal{C} \setminus \overline{\supp(D)}^{Zar}$. Then $\mathcal{U}$ is an affine model of $U$ over $\mathcal{B}$ for which all the fibers are smooth affine curves.

Every rational point  $x\in C(K) \setminus f^{-1}(\{0,1,\infty\})$ corresponds to a morphism  $\Spec(K) \longrightarrow C$ and extends to a section $\sigma_x : \mathcal{B} \longrightarrow \mathcal{C}$. The point $f(x) \in \espproj^1_{K}$ extends to a section $\sigma_{f(x)} : \mathcal{B} \longrightarrow \espproj^1_{O_{K,T}}$. Moreover $\tilde{f}(\sigma_x(\vp)) = \sigma_{f(x)}(\vp)$.

If $D=\sum n_x (x)$ is a divisor of $C$ defined over $K$ (i.e. the sum is taken over some points $x \in C(K)$), let us denote  $\tilde{D} = \sum n_x \sigma_{x}(\vp)$ the reduced divisor. (A priori, the point $\sigma_{x}(\vp) \in \mathcal{C}_{\vp}$ and, $\sigma_{x}(\vp) \in \mathcal{U}_{\vp}$ if and only if, $x$ is in $U$ and is $\vp$-integral). Moreover, $\supp(\tilde{D}) \cap \mathcal{U}_{\vp} = \supp(\widetilde{D|_U})$.

Set $H_{x}= \{\vq / \esp v_{\vq}(f(x))>0  \esp \textrm{or} \esp v_{\vq}(f(x))<0  \esp \textrm{or} \esp v_{\vq}(1-f(x))>0 \}$. Note that $\{\vq / \esp \card(\{v_{\vq}(r), v_{\vq}(1-r), v_{\vq}(1)\}) \geq 2\} = H_{x}$. Then 
$$\rad(P_{x}) = \sum_{\vq \in H_{x}} \log N(\vq).$$

We then have
\begin{displaymath}
 \begin{array}{rcl}
         
\vp \in H_{x}  &  \iff &  \sigma_{f(x)}(\vp) = \tilde{0} \esp \textrm{or} \esp \sigma_{f(x)}(\vp) = \tilde{\infty} \esp \textrm{or}  \esp \sigma_{f(x)}(\vp) = \tilde{1},\\

             &  \iff &  \tilde{f}(\sigma_{x}(\vp)) = \tilde{0} \esp \textrm{or} \esp \tilde{f}(\sigma_{x}(\vp)) = \tilde{\infty} \esp \textrm{or}  \esp \tilde{f}(\sigma_{x}(\vp)) = \tilde{1}\\

              & \iff &  \sigma_{x}(\vp) \in \supp (\tilde{D_{0}}) \cup \supp (\tilde{D_{\infty}}) \cup \supp (\tilde{D_{1}}),

 \end{array}            
 \end{displaymath}
where $\tilde{D_{0}}, \tilde{D_{1}}$ et $\tilde{D_{\infty}}$ are divisors of the special fiber $\mathcal{C}_{\vp}$.

If we suppose that $x$ is an $S$-integral point of $U$ and that $\vp\notin S$, then the reduced point is in the reduced affine curve, i.e. $\sigma_{x}(\vp) \in \mathcal{U}_{\vp}$.
In this case, we then have
\begin{displaymath}
 \begin{array}{rcl}

 \vp \in H_{x} \setminus S 
             &  \Longrightarrow &  \sigma_{x}(\vp) \in \supp (\widetilde{D_{0}|_U}) \cup  \supp (\widetilde{D_{\infty}|_U})  \cup \supp (\widetilde{D_{1}|_U})\\
             &  \iff &   \exists  \mu \in \{1, \ldots, k \}, \esp   \sigma_{x}(\vp) \in \supp (\widetilde{M_{\mu}})\\ 
             &  \iff &  \exists  \mu \in \{1, \ldots, k \}, \esp \ord_{\sigma_{x}(\vp)} (\widetilde{M_{\mu}}) \geq 1. 

 \end{array}            
 \end{displaymath}

By definition, $v_{\vp}(m_{\mu}(x)) = \frac{\ord_{\vp}(m_{\mu}(x))}{d_{\vp}} \log N(\vp)$, where $d_{\vp}$ is the local degree.
Since $\ord_{\vp}(m_{\mu}(x)) \geq \mult_{\sigma_{x}(\vp)}(\widetilde{m_{\mu}})= (2g+1) \ord_{\sigma_{x}(\vp)} (\widetilde{M_{\mu}}) \geq  2g+1$, then 
$$v_{\vp}(m_{\mu}(x)) \geq \frac{ 2g+1}{d_{\vp}} \log N(\vp) \esp \textrm{i.e.} \esp \log N(\vp)  \leq \frac{d_{\vp}}{2g+1}v_{\vp}(m_{\mu}(x)).$$

Therefore, if $x \in U(O_{K,S})$ et $\vp \in (\mathcal{B} \cap  H_{x})\setminus S$, then
 $$\log N(\vp) \leq \frac{d_{\vp}}{2g+1} \sum_{\nu = 1}^k \max \{0,v_{\vp}(m_{\nu}(x)) \}.$$
\hfill $\Box$

\noindent \textbf{Proof of the proposition \ref{radical}.}

Lemma \ref{log} gives us an upper bound for $\log N(\vp)$ for the prime ideals $\vp$ contributing to the radical of $P_{x}$ (i.e. for $\vp \in H_{x}$) for which $f(x)$ is $\vp$-integral and the curve $C$ and the functions $f$ and $m_{\nu}$ have good reduction. Denote by $\Sigma_{S}$ (respectively by $\Sigma_{T}$) the sum of $\log N(\vp)$ taken over $\vp \in S$ (respectively over $\vp \in T$). Applying lemma \ref{log} we then obtain
$$ \rad(P_{x})   \leq  \sum_{\vp\in (H \cap \mathcal{B})\setminus S} \frac{d_{\wp}}{2g+1} \sum_{\nu =1}^k \max \{0, v_{\vp}(m_{\nu}(x)) \} +  \Sigma_S  + \Sigma_T.$$
Since $(H \cap \mathcal{B})\setminus S$ is included in the set of all the valuations of $K$ and $\frac{d_v}{2g+1} \sum_{\nu = 1}^k \max \{0, v(m_{\nu}(x)) \}$ is always positive, we obtain
 $$\rad(P_{x}) \leq \sum_{\nu =1}^k \frac{1}{2g+1} \sum_{v \in M_K} d_v \max\{ 0,v(m_{\nu}(x)) \} + \Sigma_S + \Sigma_T. $$
But $\sum_{v \in M_K} d_v \max\{ 0,v(m_{\nu}(x)) \} =  h_K(m_{\nu}(x))$. Since for every $\nu \in \{1, \ldots, k\}$, $h_{K,m_{\nu}}$ and $h_{K,f}$ are two height functions on the curve $C$, by the height theory we have
$$ h_K(m_{\nu}(x)) = \frac{\deg(m_{\nu})}{\deg(f)} h_K(f(x)) + O\left(\sqrt{h_K(f(x))}\right),$$
where the implicit constant in the $O$ depends only on the degree of the number field $K$ and on the curve $C$. Since $\deg(m_{\nu}) = (2g+1)\deg(M_{\nu}) = (2g+1) d_{\nu}$ and $\deg(f)=d$, we have
$$\rad(P_{x}) \leq \frac{1}{d} \left(\sum_{\nu =1}^k d_{\nu} \right) h_K(f(x)) + c_1 \sqrt{h_K(f(x))} + \Sigma_S + \Sigma_T, $$
where the constant $c_{1}$ depends only on the degre of $K$ and on the curve $C$.   

Since for all $\nu \in \{1, \ldots,k\}$, $d_{\nu} = \deg(M_{\nu}) = \card(\supp(M_{\nu}))$, we conclude applying the lemma \ref{lemma.belyi}:  
$$\sum_{\nu =1}^k d_{\nu} = \card(\{x\in U(\Qbarre) / f(x)\in \{0, 1, \infty\}  \}) = \chi(U) + d.$$

\hfill $\Box$

\end{document}